\documentclass[a4paper,12pt]{article}
\usepackage{amssymb}\usepackage{amsmath}
\title{Normalisateurs et groupes d'Artin-Tits de type sph\'erique}
\author{Eddy Godelle}
\newbox\preuvebox
\def\preuve{\futurelet\next\lookforbracket}
\def\lookforbracket{\ifx\next[\let\go\usespecialterm
\else\let\go\relax
\ifvmode\vskip-\lastskip\fi\global\setbox\preuvebox=\vbox
\bgroup%
\noindent{\\\\\sc Preuve:}%
\enskip\relax\fi\ignorespaces\go
}
\def\usespecialterm[#1]{\ifvmode\vskip-\lastskip\fi
\global\setbox\preuvebox=\vbox\bgroup%
\noindent\hskip\parindent%
{\\\noindent\sc Preuve {\bf#1}:}\textnormal\ \ \relax
\ignorespaces}
\def\endpreuve{\hbadness10000\parfillskip=0pt\egroup
\unvbox\preuvebox
\setbox0=\lastbox
\ifdim\ht0>1pt 
\vskip-2pt
\noindent
\hbox to\textwidth{\vbox{
\parindent=0pt
\unhbox0 $\square$ \hss}\hss}
\relax%
\fi}
\newtheorem{Lemme}{Lemme}[section]

\newtheorem{Def}[Lemme]{D\'efinition}
\newtheorem{Prop}[Lemme]{Proposition}

\newtheorem{Rem}[Lemme]{Remarque}

\newtheorem{The}[Lemme]{Th\'eor\`eme}

\begin{document} \maketitle
\begin{abstract}\noindent  Let $(A_S,S)$ be an Artin-Tits and $X$ a subset of $S$ ; denote by $A_X$ the subgroup of $A_S$ generated by $X$. When $A_S$ is of spherical type, we prove that
 the normalizer and the commensurator of
 $A_X$ in $A_S$ are equal and are
 the product of $A_X$ by the quasi-centralizer of $A_X$ in $A_S$. Looking to the associated monoids  $A_S^+$ and $A_X^+$, we describe the quasi-centralizer of $A_X^+$ in $A_S^+$ thanks to results in Coxeter groups. These two results generalize earlier results of Paris (\cite{Par1}). Finaly, we compare, in the spherical case, the normalizer of a parabolic subgroup in the Artin-Tits group and in the Coxeter group.\\
2000 Mathematics Subject Classification: 20F36.\end{abstract}
\section* {Introduction}
 Soit $S$ un ensemble fini et $M= (m_{s,t})_{s,t\in S}$ une matrice
 sym\'etrique avec $m_{s,s} = 1$
pour $s\in S$ et $m_{s,t} \in \mathbb{N}\cup \{\infty\}-\{0,1\}$
pour $s\not=t$ dans $S$.  Le syst\`eme d'Artin-Tits associ\'e \`a $M$
 est la paire $(A_S,S)$  o\`u $A_S$ est le groupe d\'efini par la pr\'esentation de groupe suivante :
$$A_S = \langle S|\underbrace{sts\cdots}_{m_{s,t}\ termes} =
\underbrace{tst\cdots}_{m_{s,t}\ termes};\ \forall s,t\in S, s\not= t\
et\ m_{s,t}\not=\infty \rangle.$$ Le groupe $A_S$ est appel\'e un groupe
 d'Artin-Tits et les relations $\underbrace{sts\cdots}_{m_{s,t}\ termes} =
\underbrace{tst\cdots}_{m_{s,t}\ termes}$ sont appel\'ees les ``relations
 de tresses''.\\ Par exemple, si
$S = \{s_1,\cdots, s_n\}$ avec $m_{s_i,s_j}= 3$ pour $|i-j| =1$ et
$m_{s_i,s_j} = 2$ sinon, alors le groupe d'Artin-Tits associ\'e  est le
 groupe des tresses \`a $n+1$ brins.\\ Si l'on ajoute les relations $s^2 = 1$ \`a la pr\'esentation de
$A_S$ on obtient le groupe de Coxeter $W_S$ associ\'e \`a $A_S$ $$W_S = \langle S|\forall s\in S,\ s^2 = 1\ ; \forall s,t\in S, s\not= t\
et\ m_{s,t}\not=\infty,\  \underbrace{sts\cdots}_{m_{s,t}\ termes} =
\underbrace{tst\cdots}_{m_{s,t}\ termes} \rangle.$$  On dit
 que $A_S$, ou que $S$,
est de type sph\'erique si $W_S$ est fini.\\ Un sous-groupe $A_X$ de $A_S$ engendr\'e par une partie $X$
 de $S$ est appel\'e un sous-groupe parabolique standard; un
 sous-groupe conjugu\'e \`a un tel sous-groupe est appel\'e
un sous-groupe parabolique. Van Der Lek a montr\'e dans \cite{VdL} que
 $(A_X,X)$ est canoniquement 
isomorphe au syst\`eme d'Artin-Tits associ\'e \`a la matrice $(m_{s,t})_{s,t\in
X}$. Pour une telle partie $X$ de $S$, on appelle centralisateur, quasi-centralisateur,
normalisateur et commensurateur de $A_X$ dans $A_S$ les sous-groupes de $A_S$ respectifs $Z_{A_S}(A_X) = \{g\in A_S| \forall s\in X,\  gs = sg\}$, $QZ_{A_S}(A_X) = \{g\in A_S| gX = Xg\}$, 
$N_{A_S}(A_X) = \{g\in A_S| gX \subset A_Xg\}$, $Com_{A_S}(A_X)
  = \{g\in A_S| gA_Xg^{-1}\cap A_X$ est d'indice fini dans $A_X$ et dans
$gA_Xg^{-1}$ $\}$.\\ 

\noindent Les principaux r\'esultats relatifs aux normalisateurs sont les suivants :
\begin{The} \label{THdeb2}
Soit $(A_S,S)$ un syst\`eme d'Artin-Tits de type sph\'erique et $X\subset~S$. On choisit $\epsilon\in \{1,2\}$ minimal tel que $\Delta_X^\epsilon$ est central dans $A_X$ ({\it cf.} lemme \ref{lemclassi}(iii)). Alors  $$Z_{A_S}(\Delta_X^\epsilon) =  Com_{A_S}(A_X) = N_{A_S}(A_X) = A_X \cdot QZ_{A_S}(A_X).$$
\end{The}
Ce th\'eor\`eme a \'et\'e d\'emontr\'e par Paris ( \cite{Par1}
th\'eor\`eme 5.1) avec l'hypoth\`ese suppl\'ementaire que $X$ est
ind\'ecomposable ({\it cf.} la d\'efinition suivant le th\'eor\`eme \ref{THcox}). \begin{The}\label{maconj}
									       Soit $(A_S,S)$ un syst\`eme d'Artin-Tits de type sph\'erique, soit $G = gA_Xg^{-1}$ un sous-groupe parabolique de $A_S$ et $Y\subset~S$. Si $G\subset A_Y$ alors $G$ est un sous-groupe parabolique de $A_Y$: $$gA_Xg^{-1}\subset A_Y \Rightarrow \exists R\subset Y,\exists y\in A_Y \textrm{ tels que } gA_Xg^{-1} = yA_Ry^{-1}.$$  \end{The}
\begin{The}\label{THcox}Soit $(A_S,S)$ un syst\`eme d'Artin-Tits de type sph\'erique et $X$ une partie de $S$, alors  $$N_{A_S}(A_X)/(Z_{A_S}(A_X)\cdot A_X)\simeq N_{W_S}(W_X)/(Z_{W_S}(W_X)\cdot W_X).$$
\end{The}

\noindent Soit $(A_S,S)$ un syst\`eme d'Artin-Tits (quelconque). On dit que $A_S$ (ou simplement $S$) est d\'ecomposable s'il existe une partition non triviale de $S = S_1\cup S_2$ telle que $A_S$ est canoniquement isomorphe au
 produit direct $A_{S_1}\times A_{S_2}$ ; c'est \`a dire que $\forall
 s\in S_1,\forall t\in S_2$, $m_{s,t} = 2$. Si $S$ n'est pas
 d\'ecomposable, on dit qu'il est ind\'ecomposable. Les composantes ind\'ecomposables de $S$ sont les
 sous-ensembles ind\'ecomposables maximaux de $S$.\\
 Soit $A^+_S$ le mono\"\i de engendr\'e par $S$ dans $A_S$. Alors ({\it cf.} \cite{BrS}, \cite{Par4}) $A^+_S$ poss\`ede la pr\'esentation de mono\"\i de suivante: $$A^+_S = \langle S|\underbrace{sts\cdots}_{m_{s,t}\ termes} =
\underbrace{tst\cdots}_{m_{s,t}\ termes};\ \forall s,t\in S, s\not= t\
et\ m_{s,t}\not=\infty \rangle^+.$$ Si $T,X,Y \subset S$, on appelle conjugateur dans $A_T$ de $X$ en $Y$  et conjugateur positif dans $A_T^+$ de $X$ en $Y$ les ensembles respectifs $Conj(T;X,Y) = \{g\in A_T | gX = Yg\}$ et $Conj^+(T;X,Y) = \{g\in A^+_T | gX = Yg\}$. En particulier, $Conj(S;X,X) = QZ_{A_S}(A_X)$.    
\begin{Def}\label{categ+}
Soit $X \subset S$ et $t\in S $ tels que la composante ind\'ecomposable $X(t)$ de $X\cup \{t\}$ contenant $t$ est de type sph\'erique. 
Si $t\not\in X$, on pose $d_{X,t} = \Delta_{X(t)}\Delta_{X(t)-\{t\}}^{-1}$ ; sinon, on pose $d_{X,t} = \Delta_{X(t)}$. Dans les deux cas, il existe une unique partie $Y$ de $X\cup \{t\}$ telle que $Yd_{X,t} = d_{X,t} X$. On dira alors que $d_{X,t}$ est un $Y$-ruban-$X$ \'el\'ementaire positif.\\
Si $X,Y\subset S$, on dira que $g\in A^+_S$ est un $Y$-ruban-$X$ positif si $g =  g_n\cdots g_1$ o\`u $g_i$ est un $X_i$-ruban-$X_{i-1}$ positif \'el\'ementaire, $X_0 = X$ et $X_n = Y$. \end{Def}
On la r\'esultat suivant: 
\begin{The}\label{cordeb1}
Soit $(A_S,S)$ un syst\`eme d'Artin-Tits et $X,Y\subset S$, alors
$$\displaylines{g\in Conj^+(S;X,Y)\iff g= g_n\cdots g_1\textrm{ o\`u }g_i\textrm{
est un }X_i\textrm{-ruban-}X_{i-1}\hfill\cr\hfill\textrm{ positif \'el\'ementaire, }X_0 = X\textrm{ et }X_n = Y.}$$\end{The}
C'est \`a dire que les \'el\'ement de $Conj^+(S;X,Y)$ sont exactement les $Y$-ruban-$X$ positifs. Ce r\'esultat a \'et\'e prouv\'e par Paris dans \cite{Par1} pour les groupes
d'Artin-Tits de type sph\'erique.\\
 
\noindent Dans la premi\`ere partie nous rappelons les r\'esultats relatifs aux groupes d'Artin-Tits dont nous aurons besoin, dans la seconde partie, nous prouvons les th\'eor\`emes \ref{THdeb2} et \ref{maconj} gr\^ace \`a la proposition clef \ref{Propclef}. Dans la troisi\`eme partie, on \'etablit le th\'eor\`eme \ref{cordeb1} apr\`es un d\'etour par les groupes de Coxeter. Enfin nous terminons en prouvant le th\'eor\`eme \ref{THcox} dans la derni\`ere partie.
\section{G\'en\'eralit\'es}
Soit $(A_S,S)$ un syst\`eme d'Artin-Tits et $A^+_S$ le mono\"\i de engendr\'e par $S$ dans $A_S$. Alors ({\it cf.} \cite{BrS}, \cite{Par4}) $A^+_S$ poss\`ede la pr\'esentation de mono\"\i de suivante: $$A^+_S = \langle S|\underbrace{sts\cdots}_{m_{s,t}\ termes} =
\underbrace{tst\cdots}_{m_{s,t}\ termes};\ \forall s,t\in S, s\not= t\
et\ m_{s,t}\not=\infty \rangle^+.$$
Puisque les relations de la pr\'esentation de $A^+_S$ sont homog\`enes, on peut munir $A^+_S$ d'un unique morphisme de mono\"\i des, appel\'e longueur et not\'e $\ell : A^+_S\to \mathbb{N}$, v\'erifiant $\ell(s) = 1$ pour $s\in S$.\\
On d\'esignera par $\prec$ et $\succ$ respectivement la division \`a droite et la division \`a gauche dans $A^+_S$. On note $a\land_\prec b$ et $a\lor_\prec b$ ({\it resp.} $a\land_\succ b$ et $a\lor_\succ b$) les pgcd et ppcm pour $\prec$ ({\it resp.} $\succ$) de $a$ et $b$ dans $A^+_S$ lorsque ceux-ci existent. 
\begin{Lemme}[\cite{BrS}, \cite{Mic} proposition 2.4 et 2.6] \label{lemmic} Soit $(A_S,S)$ un syst\`eme d'Artin-Tits.\\
(i) $A_S^+$ est  simplifiable.\\
(ii) toute partie finie de $A_S^+$ poss\`ede un pgcd pour $\prec$(et pour
 $\succ$).\\
(iii) une partie finie de $A_S^+$ poss\`ede un ppcm pour $\prec$ (\it{resp.} pour $\succ$) si et seulement si elle poss\`ede un multiple commun pour $\prec$  (\it{resp.} pour $\succ$).
\end{Lemme}   
Si $X,Y\subset S$, on dira que $g\in A^+_S$ est $X$-r\'eduit ({\it resp.} r\'eduit-$Y$) s'il n'est divisible pour $\prec$ ({\it resp.} $\succ$) par aucun \'el\'ement de $X$ ({\it resp.} $Y$). Enfin, on dira que $g$ est $X$-r\'eduit-$Y$ s'il est \`a la fois $X$-r\'eduit et r\'eduit-$Y$.\\ 

Lorsque $S$ est de type sph\'erique, le mono\"\i de $A^+_S$ est un mono\"\i de de Garside (\cite{DeP}). En particulier, l'ensemble $S$ poss\`ede alors un ppcm que l'on note $\Delta_S$ ;  on a \'egalement:  
\begin{Lemme} \label{lemclassi}
Soit $A_S$ un groupe d'Artin-Tits de type sph\'erique.\\ (i) {\bf (\cite{Del} paragraphe 4)} tout \'el\'ement $g$ de $A_S$ s'\'ecrit $g = g_1\Delta_S^n$ avec $g_1\in A_S^+$, et $n$ dans $\mathbb{Z}$.\\
(ii) {\bf(\cite{Cha1} th\'eor\`eme 2.6 et $\cite{Cha2}$ lemme 4.4)} Si $g\in A_S$ alors il existe $a$ et $b$ dans $A_S^+$ uniques tels que $g = ab^{-1}$ et $a\land_\succ b = 1$. De plus si $c\in A_S^+$ est tel que $gc\in A_S^+$ alors $b\prec c$.\\
(iii) La conjugaison par $\Delta_S$ est un automorphisme de $A_S$ d'ordre 1 ou 2 qui stabilise $S$, en particulier il induit un automorphisme de $A^+_S$. \\  
\end{Lemme} 
Nous appellerons  la d\'ecomposition $g = ab^{-1}$ l'\'ecriture
normale \`a droite de $g$, on peut de m\^eme d\'efinir l'\'ecriture
normale \`a gauche.\\ On notera que si $S$ est de type sph\'erique et $X\subset S$ alors $X$ est aussi de type sph\'erique ; dans ce cas si $g\in A_X$ et $g = ab^{-1}$ est l'\'ecriture
normale \`a droite de $g$ dans $A_S$ alors par l'unicit\'e d'une telle \'ecriture c'est aussi son \'ecriture
normale \`a droite dans $A_X$ (en particulier $a,b\in A_X^+$).
\begin{Rem} Une r\'e\'ecriture tr\`es utile du (ii) du lemme \ref{lemclassi} est la suivante~: sous les hypoth\`eses de ce lemme, si $g = ab^{-1} = uv^{-1}$ avec $a,b,u,v\in A^+_S$ et $a\land_\succ b = 1$ alors il existe $\alpha\in A^+_S$ tel que $u = a\alpha$ et $v = b\alpha$. 
\end{Rem}
Soit $A_S$ un groupe d'Artin-Tits. On note $p^+: A_S^+\to W_S$ la surjection canonique. On d\'efinit la longueur $\ell(z)$ d'un \'el\'ement $z\in W_S$ comme le minimum des longueurs des \'el\'ements de $(p^+)^{-1}(z)$; on a donc $\ell(g) \geq \ell(p^+(g))$ pour tout $g\in A^+_S$.\\
Il est connu (\textit{cf.} \cite{Mic}  sect. 1) que $p^+$
 donne lieu -- gr\^ace au lemme d'\'echange dans les groupes de Coxeter -- \`a une section $\pi$ dont l'image, not\'ee $A_{S,red}$, est
 form\'ee des \'el\'ements de $A_S^+$ qui ont m\^eme longueur que leur
 image dans $W_S$ par $p^+$. Cet ensemble est en particulier stable par division \`a gauche et la division \`a droite.
Les \'el\'ements de $A_{S,red}$\index{$A_{S,red}$} seront dit ``r\'eduits''.
\begin{Lemme}[\cite{BrS} lemme 3.4] \label{lemred2}Soit $(A_S,S)$ un
 syst\`eme d'Artin-Tits, soit
 \linebreak $x\in A_{S,red}$ et $s\in S$; si $sx\not\in A_{S,red}$ alors $s\prec x$.
\end{Lemme}
\begin{Lemme} [(\cite{Thu}, \cite{Mic} proposition 1.7 et 2.1)]  \label{lemmic2} Il existe une unique\index{$\alpha(g)$}
application $\alpha~: A_S^+ \rightarrow A_{S,red}$ qui est l'identit\'e sur
$A_{S,red}$, v\'erifie $\alpha(gh) = \alpha(g\alpha(h))$
pour $g,h \in A_S^+$ et $\alpha(ab) = ac$ pour $a,b\in A_{S,red}$ avec $c$ le plus grand \'el\'ement de l'ensemble $\{d\in A_{S,red}| d\prec b\ et\ ad\in A_{S,red}\}$. De plus $\alpha(g)$ est le plus grand \'el\'ement de l'ensemble $\{c\in A_{S,red}| c\prec g\}$ o\`u $\prec$ d\'esigne la division \`a gauche.\end{Lemme}
\begin{Prop} \label{Prfornor}Soit $A_S^+$ un mono\"\i de d'Artin-Tits et $g\in A^+_S-\{1\}$ ; il existe une unique suite finie $(g_1,\cdots, g_n)$ d'\'el\'ements de $A_{S,red}-\{1\}$ telle que $g = g_1g_2\cdots g_n$ et $g_i = \alpha(g_i\cdots g_n)$ pour tout $i\in \{1,\cdots, n\}$. \end{Prop}
On dira que cette suite est la d\'ecomposition normale d'Adyan ({\it cf.} \cite{Ady}) de $g$ (\`a gauche).\\ 

Dans la derni\`ere partie nous aurons besoin de la notion de cha\^\i ne (\`a gauche): 
\begin{Def} \label{defprcha3}Soit $A^+_S$ un mono\"\i de d'Artin-Tits. Soit $s,t\in S$ distincts et $u\in \{s,t\}$. Soit $C\in A^+_S$ ;\\ (i) on dit que $C$ est une $s$-cha\^\i ne-$u$ simple, si $C = 1$ et $s = u$ ou bien si $C = \underbrace{tst\cdots}_{k\ termes}$ avec $k < m_{s,t}$ et $Cu = \underbrace{tst\cdots}_{k+1\ termes}$.\\(ii) On dit que $C$ est une $s$-cha\^\i ne-$t$ s'il existe une suite $s_0 = s,s_1,\cdots, s_k = t$ d'\'el\'ements de $S$ tels que $C =  C_1\cdots C_k$ o\`u $C_i$ une $s_{i-1}$-cha\^\i ne-$s_i$ simple pour $i\in\{1,\cdots, k\}$.\\ Dans ce cas, on dit que $s$ est l'origine de $C$ et $t$ son but.  
\end{Def}
On parle parfois  simplement de $s$-cha\^\i ne ou de cha\^\i ne-$t$ pour parler d'une $s$-cha\^\i ne-$t$ dont on ne pr\'ecise pas l'une des extr\'emit\'es. 
\begin{Prop}[\cite{BrS} lemmes 3.1 et 3.2] \label{Prprcha2}Soit $A^+_S$ un mono\"\i de d'Artin-Tits. Soit $s\in S$ et $g,h\in A_S^+$; alors\\ (i) si $g$ est une $s$-cha\^\i ne-$t$ et $s\prec gh$ alors $t\prec h$.\\ (ii) $s$ ne divise pas $h$ \`a gauche  si et seulement si $h = gk$ o\`u $g$ est une $s$-cha\^\i ne-$t$ pour un certain $t\in S$ et, ou bien $k =1$  ou bien $k= rk_1$ avec $r\in S$ et $k_1\in G^+$ tels que  $r\lor_\prec t$ n'existe pas.        
\end{Prop} 
Par sym\'etrie des relations, on peut \'enonc\'e un r\'esultat similaire pour $\succ$.
\section{Preuve des th\'eor\`emes \ref{THdeb2} et \ref{maconj}} 
Les d\'emonstrations des th\'eor\`emes $\ref{THdeb2}$ et \ref{maconj} reposent essentiellement sur la proposition suivante :\begin{Prop}[proposition clef]\label{Propclef} Soit $(A_S,S)$ un syst\`eme d'Artin-Tits de type sph\'erique et $X,Y \subset S$. Soit $k\in \mathbb{Z}-\{0\}$ et $g\in A_S$. Les assertions suivantes sont \'equivalentes.\\
(1) $gA_X g^{-1}\subset A_Y$;\\
(2) $g\Delta_X^kg^{-1} \in A_Y$;\\
(3) $g = yx$ avec $y\in A_Y$, $x\in Conj(S;X,R)$ pour une partie $R\subset Y$.\\    
\end{Prop}
Il est clair que l'implication $(1)\Rightarrow (3)$ de cette proposition entra\^\i ne le th\'eor\`eme \ref{maconj}.  Afin d'\'etablir la proposition \ref{Propclef}, nous commen\c cons par prouver le lemme suivant. 
\begin{Lemme}\label{lemcentdel} Soit $(A_S,S)$ un groupe d'Artin-Tits; soit $X,Y\subset S$
de type sph\'erique et $g\in A_S^+$ r\'eduit-$X$; Alors~: $$\left(\exists k\in \mathbb{N}^*\textrm{ tel que } g\Delta_X^kg^{-1}\in A^+_Y \right)\Rightarrow g\in Conj^+(S;X,R) \textrm{ pour un } R\subset Y.$$  
\end{Lemme}
\begin{preuve} On montre le r\'esultat par r\'ecurrence sur $\ell(g)$. Si
 $\ell(g) = 0$ alors $g = 1$ et le r\'esultat est vrai avec $R= X$. Supposons donc
 $\ell(g)\geq 1$ et que $\forall X',Y'\subset S$ de type sph\'erique et $\forall g'\in A_S^+$ r\'eduit-$X'$ \newpage $$\displaylines{(\ell(g')<\ell(g)\textrm{ et }\exists j\in \mathbb{N}^*\textrm{ avec }  g'\Delta_{X'}^j{g'}^{-1}\in A^+_{Y'}) \Rightarrow \hfill\cr\hfill g'\in Conj^+(S;X',R') \textrm{ pour un } R'\subset Y'.}$$ Notons $z = g\Delta_X^kg^{-1}\in A_Y^+ $.  Puisque $g$ est r\'eduit-$X$, il existe
 $t\in S-X$ tel que $g\succ t$. L'\'egalit\'e
 $zg = g\Delta_{X}^k$ montre que
 $zg$ est divisible (\`a droite) par tous les \'el\'ements de
 $X\cup\{t\}$. Donc $X\cup\{t\}$ est de type sph\'erique et
 $zg\succ\Delta_{X\cup\{t\}}$. Notons $d_{X,t} = \Delta_{X\cup\{t\}}\Delta_X^{-1}$. On a
 $g\Delta^{k-1}_X \succ d_{X,t}$ et comme
 $\Delta_{X\cup\{t\}}\Delta_X^{-1}\succ t$, on trouve de nouveau que $g\Delta_X^{k-1}\succ \Delta_{X\cup\{t\}}$. Par une r\'ecurrence imm\'ediate sur $k$  on trouve que $g \succ d_{X,t}$~: on a $g = g' d_{X,t}$ avec $g'\in A_S^+$. D'o\`u $zg'd_{X,t} = g'd_{X,t}\Delta_{X}^k$ et par simplifiabilit\'e  $zg' = g'\Delta_{R}^k$ o\`u $R\subset S$ est tel que $Rd_{X,t} = d_{X,t}X$. Maintenant l'hypoth\`ese de r\'ecurrence appliqu\'ee \`a $g'$ qui est bien r\'eduit-$R$, implique que $g'$ est dans $Conj^+(S;R,R')$ avec $R'\subset Y$. Enfin, comme $d_{X,t}\in Conj^+(S;X,R)$ (\textit{cf.} lemme 1.2(iii)), on trouve que $g\in Conj^+(S;X,R')$ pour $R'\subset Y$. 
\end{preuve}
Comme nous allons le voir dans le paragraphe suivant, les \'el\'ements de la forme $d_{X,t}$ sont en fait ``les briques \'el\'ementaires '' pour construire les rubans qui servent \`a d\'ecrire $Conj^+(S;X,Y)$.\\
\begin{preuve}[de la proposition \ref{Propclef}]
Il est clair que $(3)\Rightarrow (1)\Rightarrow (2)$. Nous devons donc juste prouver que $(2)\Rightarrow (3)$. Soit $g\in A_S$ tel que $g\Delta_X^kg^{-1} \in A_Y$ ; on peut sans restriction supposer que $k\in\mathbb{N}^*$ puis que $k$ est pair. Par le lemme \ref{lemclassi}(i), Il existe $n\in \mathbb{N}$ tel que $g = h\Delta_S^{2n}$ avec $h\in A_S^+$ (quitte \`a remplacer dans ce lemme $g_1$ par $g_1\Delta_S$, on peut toujours supposer que la puissance de $\Delta_S$ est paire). Puisque $\Delta_S^2\in Z_{A_S}(A_X)$ (lemme \ref{lemclassi}(iii)), on a $h\Delta_X^kh^{-1} \in A_Y$ et si on \'ecrit $h=abc$ avec $c\in A_X^+$, $a\in A_Y^+$ et o\`u $b$ est $Y$-r\'eduit-$X$ (cette \'ecriture n'est pas unique, mais il nous suffit d'en choisir une). Puisque $k$ est pair, on a aussi $b\Delta_X^k b^{-1} \in A_Y$. Par cons\'equent, l'\'ecriture normale de $b\Delta_X^ kb^{-1}$ dans $A_S$
 est de la forme $uv^{-1}$  avec $u,v$ dans $A_Y^+$ et $u\land_\succ v = 1$. En
 appliquant le lemme \ref{lemclassi}(ii) \`a l'\'egalit\'e $b\Delta_X^ kb^{-1} = uv^{-1}$, on trouve que $b = v\alpha$ avec  $\alpha$ dans $A_S^+$. Mais $b$ est $Y$-r\'eduit, donc $v = 1$ et $b\Delta_X^ kb^{-1}\in A^+_Y$. Par le lemme \ref{lemcentdel}, cela
 implique que $b\in Conj^+(S;X,R)\subset Conj(S;X,R)$ pour $R\subset Y$. 
Finalement on trouve $g = yx$ avec $y = abcb^{-1}\in A_Y$ et $x = b\Delta^{2n}\in Conj(S;X,R)$ pour un $R\subset Y$.\end{preuve}
\begin{preuve}[du th\'eor\`eme \ref{THdeb2}] Il est clair que $A_X \cdot QZ_{A_S}(A_X)\subset N_{A_S}(A_X)\subset Com_{A_S}(A_X)$ et que $A_X \cdot QZ_{A_S}(A_X)\subset Z_{A_S}(\Delta_X^\epsilon)$. Il suffit donc de montrer $Com_{A_S}(A_X)\subset A_X \cdot QZ_{A_S}(A_X)$ et $Z_{A_S}(\Delta_X^\epsilon)\subset A_X \cdot Q_{A_S}(A_X)$. Par la proposition \ref{Propclef}, on a l'inclusion $$\{g\in A_S; \exists k\in \mathbb{N}^*\textrm{ tel que } g\Delta_X^kg^{-1}\in X\} \subset A_X\cdot QZ_{A_S}(A_X) = QZ_{A_S}(A_X) \cdot A_X.$$ On a alors imm\'ediatement que $Z_{A_S}(\Delta_X^\epsilon)\subset A_X\cdot QZ_{A_S}(A_X)$. D'autre part si $g\in Com_{A_S}(A_X)$, alors pour tout $z\in A_X$, il existe $k\in \mathbb{N}^*$ tel que $gz^kg^{-1}\in A_X$. Il suffit de choisir $z = \Delta_X$ pour terminer la preuve. 
\end{preuve}
\section{Conjugateur dans le mono\"\i de}
\noindent Commen\c cons par regarder le cas des groupes de Coxeter. Soit $(W_S,S)$ un syst\`eme de Coxeter. On appelle base de racines\index{racine} de $W_S$ un triplet $(E;(.,.);\Pi)$ o\`u $E$ est un espace vectoriel de dimension $|S|$, o\`u  $(.,.)$ est une forme bilin\'eaire sym\'etrique sur $E$ et $\Pi=\{e_s; s\in S\}$ est une base de $E$ telle que $(e_s,e_t) = -cos(\pi/m_{s,t})$ si $m_{s,t}$ est fini et $(e_s,e_t)\leq -1$ si $m_{s,t} = \infty$. $W_S$ agit fid\`element sur $E$ par $$s.v = v -2(v,e_s)e_s\textrm{ pour $s\in S$ et }v\in E.$$ L'ensemble $\Phi = \{w\cdot e_s | s\in S$ et $w\in W_S \}$ est le syst\`eme de racines de $W_S$ dans $E$ et les \'el\'ements de $\Pi$ s'appellent les racines simples. Une racine est dite positive (\textit{resp.} n\'egative) si elle s'\'ecrit $\sum_{s\in S}{\lambda_s e_s}$ avec $\lambda_s \geq 0$ ({\it resp.} $\lambda_s \leq 0$) ; L'ensemble des racines positives ({\it resp.} n\'egatives) est not\'e $\Phi^+$ ({\it resp.} $\Phi^-$). On a $\Phi = \Phi^+\cup \Phi^-$ et $-\Phi^+ = \Phi^-$. Pour chaque $\alpha\in \Phi$, il existe une r\'eflexion $s_\alpha$ de $W_S$, c'est \`a dire le conjugu\'e $s_\alpha = wsw^{-1}$ d'un \'el\'ement de $S$, tel que $s_\alpha.v = v -2(v,\alpha)\alpha$ pour tout $v\in E$ ; de plus, $\alpha = w\cdot e_s$. Soit $\Theta\subset \Phi$ et notons $W_\Theta$ le sous-groupe de $W_S$ engendr\'e par $\{s_\theta ;\theta\in \Theta\}$ et $\Phi_\Theta = \{w\cdot \theta ; \theta\in \Theta$ et $w\in W_\Theta\}$. Alors $W_\Theta$ est un groupe de Coxeter dont $\Phi_\Theta$ est un syst\`eme de racines (\cite{Deo}).\\ Soit $X$ une partie de $S$.  On pose $\Pi_X= \{e_s; s\in X\}$ et on note $\omega_X$ l'\'el\'ement de plus grande longueur de $W_X$ lorsqu'il existe ; on a $p^+(\Delta_X) = \omega_X$ o\`u $p^+ : A^+_S\to W_S$ est la surjection canonique.\\
Rappelons que si $X\subset S$ et $s\in S$, on d\'esigne par $X(s)$ la composante ind\'ecomposable de $X\cup\{s\}$ qui contient $s$. Lorsque $s\in S-X$ est tel que $X(s)$ est de type sph\'erique, on pose $\nu(X,s) = \omega_{X(s)-\{s\}}\omega_{X(s)}$.  
\begin{Prop}[\cite{Deo} proposition 5.5 ; voir aussi \cite{How} lemme 5]\label{propdeo}\ \\ Soit $(W_S,S)$ un syst\`eme de Coxeter.\\
(i) Soit $X\subset S$ et $s\in S-X$ tel que $\nu(X,s)$ est d\'efini. Alors il existe $t\in X(s)$ tel que $\nu(X,s)^{-1}\Pi_X = \Pi_Y$ avec $Y = X\cup\{s\}-\{t\}$ et $t\in X(s)$. On dit que $X$ est l'origine de $\nu(X,s)$ et $Y$ son but.\\
(ii) Soit $X,Y\subset S$ et $w\in W_S$ tels que $w^{-1}\Pi_X = \Pi_Y$. Alors il existe une suite $s_0,\ldots,s_n$ d'\'el\'ements de $S$ et une suite $X_0 =X,X_1,\ldots,X_{n-1},X_n = Y$ de parties de $S$ telles que $w = \nu(X_0,s_0)\cdots \nu(X_n,s_n)$ o\`u pour $i\in\{0,\ldots,n-1\}$, le but de $\nu(X_i,s_i)$ est $X_{i+1}$ et telles que $l(w) = \sum_{i=0}^{n-1}{\ell(\nu(X_i,s_i))}$.   
\end{Prop}
\begin{Prop}[\cite{Kra} proposition 3.1.9]\label{propkracox}
Soit $(W_S,S)$ un syst\`eme de Coxeter et $X$ une partie de $S$. Alors on a $N_{W_S}(W_X) = G_X\ltimes W_X$ o\`u $G_X$ est le groupe $\{w|$ $w\Pi_X = \Pi_X\}$. 
\end{Prop}
La proposition suivante est l'argument essentiel permettant de prouver le th\'eor\`eme \ref{cordeb1}.
\begin{Prop}
Soit $A_S$ un groupe d'Artin-Tits. Soit $s,t \in S$, soit $g$ un \'el\'ement de $A_S^+$ et $g =g_1\cdots g_n$ sa d\'ecomposition normale d'Adyan. Si $gs = tg$, alors il existe une suite $s_0 =t,s_1,\cdots,s_n =s$ de $S$ telle que pour $i\in \{1,\cdots,n\}$, on a $ g_is_i = s_{i-1}g_i$.
\end{Prop} 
\begin{preuve} On fait une r\'ecurrence sur $n$. Si $n = 1$, il n'y a rien \`a montrer ; supposons donc $n\geq 2$ (en particulier $g\not= 1$). L'ensemble $Conj_{red}(g,t) = \{h\in A_{S,red} | h\prec g$ et $\exists v\in S,\ th = hv\}$ n'est pas vide : si $u\in S$ avec $u\prec g$ alors $m_{t,u}\not=\infty$ car $u$ et $t$ ont un multiple commun et $\underbrace{utu\cdots}_{m_{t,u}-1\ termes } \in Conj_{red}(g,t)$. Soit $z_1\in Conj_{red}(g,t)$ de longueur maximale ; puisque $z_1$ est r\'eduit, on a par le lemme \ref{lemmic2} que $g_1 = z_1z_2$ avec $z_2$ r\'eduit. Supposons $z_2\neq 1$ ; par d\'efinition, $tz_1 = z_1s_1$ pour un certain $s_1\in S$, et par simplifiabilit\'e dans l'\'egalit\'e $tg = gs$, on obtient $s_1z_2g_2\cdots g_n = z_2g_2\cdots g_n s$. Soit $u\in S$ avec $u\prec z_2$ et notons $z_2 = uz'_2$ avec $z'_2\in A^+_S$. De l'\'egalit\'e $s_1z_2g_2\cdots g_n = z_2g_2\cdots g_n s$ on d\'eduit que $s_1$ et $u$ ont un multiple commun ; donc $m_{s_1,u}$ est fini et par simplifiabilit\'e,  $\underbrace{us_1u\cdots }_{m_{u,s_1}-1} \prec z_2g_2\cdots g_n$. Par maximalit\'e de la longueur de $z_1$ dans $Conj_{red}(g,t)$, l'\'el\'ement  $z_1\underbrace{us_1u\cdots}_{m_{u,s_1}-1}$ ne peut \^etre r\'eduit puisqu'il divise $g$ et v\'erifie $t z_1\underbrace{us_1u\cdots}_{m_{u,s_1}-1}= z_1\underbrace{us_1u\cdots}_{m_{u,s_1}-1}v $ pour un certain $v\in S$. On en d\'eduit que $z_1\succ us_1$ par le lemme \ref{lemred2}. Notons $z'_1\in A^+_S$ tel que $z_1 = z'_1us_1$. De l'\'egalit\'e $tz_1 = z_1s_1$ on obtient apr\`es simplification $tz'_1u = z'_1us_1 = z_1$ ; ce qui implique que $z_1 = z''_1\underbrace{us_1u\cdots}_{m_{s_1,u}\ termes}$. En regroupant toutes ces \'ecritures on trouve $g_1 = z''_1\underbrace{us_1u\cdots}_{m_{s_1,u}\ termes}uz'_2$, ce qui contredit le fait que $g_1$ est r\'eduit. Donc $z_2 = 1$ et $tg_1 = g_1s_1$ ; on conclut maintenant en appliquant l'hypoth\`ese de r\'ecurrence \`a $g_2\cdots g_n$ et $s_1g_2\cdots g_n = g_2\cdots g_n s$.
\end{preuve}
\begin{Lemme}\label{lemsph1}
Soit $(A_S,S)$ un syst\`eme d'Artin-Tits de type sph\'erique. Soit $X \subset S$ ind\'ecomposable, $t\in X$ et $g$ un $Y$-ruban-$X$ positif pour un certain $Y\subset S$. Si $g\succ t$ alors $g\succ\Delta_X$.\end{Lemme}
\begin{preuve} si $g\succ t$ et $s\in S$ tel que $m_{s,t}\neq 2$, une preuve analogue \`a celle du lemme 3.3 de \cite{Mic} montre que $g\succ s$. On en d\'eduit comme dans \cite{Mic} que $g\succ \Delta_X$. \end{preuve} 
\begin{preuve}[du th\'eor\`eme \ref{cordeb1}]
Le sens $\Leftarrow$ est clair ; supposons donc $g\in Conj^+(S,X,Y)$ et notons $g = g_1\cdots g_n$ sa d\'ecomposition normale d'Adyan. Par la proposition pr\'ec\'edente, pour tout $s\in X$, il existe une suite $s_0 = s,s_1,\cdots,s_n \in Y$ de $S$ telle que pour $i\in \{1,\cdots,n\}$, $g_is_i = s_{i-1}g_i$. Il existe donc une suite $X_0 = X,X_1,\cdots,X_n = Y$ de parties de $S$ telles que pour $i\in \{1,\cdots,n\}$, $g_i\in Conj^+(S;X_{i-1},X_i)$. Il suffit donc de montrer l'implication ``$\Rightarrow$'' lorsque $g\in A_{S,red}$.  Supposons donc $g$ r\'eduit. On a alors $g = \pi(p^+(g))$ o\`u $\pi$ est la section de $p^+$ d\'efinie dans la premi\`ere partie ; on a aussi $Yp^+(g) = p^+(g)X$. Supposons pour commencer que $p^+(g)$ est r\'eduit-$X$. Soit $s\in Y$ et $t\in X$ tels que $p^+(g)t = sp^+(g)$. On a alors $t\cdot(g^{-1}\cdot e_s) = -g^{-1}\cdot e_s$. Donc $g^{-1}(e_s) \in \{e_t,-e_t\}$. Mais puisque $g$ est r\'eduit-$X$, il envoie une racine positive sur une racine positive et $g^{-1}(e_s)  = e_t $. D'o\`u $g^{-1}\Pi_Y = \Pi_X$. Par la proposition \ref{propdeo}(ii), il existe une suite $s_0,\ldots,s_n$ d'\'el\'ements de $S$ et une suite $X_0 =Y,X_1,\ldots,X_{n-1},X_n = X$ de parties de $S$ telles que $p^+(g) = \nu(X_0,s_0)\cdots \nu(X_n,s_n)$ o\`u pour $i\in\{0,\ldots,n-1\}$, le but de $\nu(X_i,s_i)$ est $X_{i+1}$ et telles que $\ell(p^+(g)) = \sum_{i=0}^{n-1}{l(\nu(X_i,s_i))}$. Puisque $\ell(p^+(g)) = \sum_{i=0}^{n-1}{l(\nu(X_i,s_i))}$, on a $g =  \pi(p^+(g)) = \pi(\nu(X_0,s_0))\cdots \pi(\nu(X_n,s_n))$. Mais $\pi(\nu(X_i,s_i)) = \Delta_{X_i(s_i)-\{s_i\}}^{-1}\Delta_{X_i(s_i)} = \Delta_{X_{i+1}(t_i)} \Delta_{X_{i+1}(t_i)-\{t_i\}}^{-1}$ o\`u $t_i$ est d\'efini par $X_i\cup\{s_i\} = X_{i+1}\cup \{t_{i}\}$. Donc $g$ est un $Y$-ruban-$X$ positif.\\
Supposons maintenant que $g$ n'est pas r\'eduit-$X$ et choisissons $t\in X$ tel que $g\succ t$. Alors par le lemme \ref{lemsph1}, $g = g_1\Delta_{X(t)}$ avec $g_1\in A^+_S$. Puisque $\Delta_{X(t)}$ est par d\'efinition un $X$-ruban-$X$ positif, on a $g_1\in Conj^+(S;X,Y)$. Par r\'ecurrence sur $\ell(g)$, on peut ainsi \'ecrire $g = g'g''$ avec $g''$ un $X$-ruban-$X$ positif et o\`u $g'\in Conj^+(S;X,Y)$  est r\'eduit (il divise $g$) et r\'eduit-$X$.  Par la premi\`ere partie de la preuve, $g'$ est $Y$-ruban-$X$. Enfin puisque $g''$ est un $X$-ruban-$X$ positif, $g = g'g''$ est un $Y$-ruban-$X$ positif.\end{preuve}
\section{Preuve du th\'eor\`eme \ref{THcox}}

Commen\c cons par quelques lemmes techniques.
\begin {Lemme}\label{corZen} Soit $(A_S,S)$ un syst\`eme d'Artin-Tits de type sph\'erique et $s,t\in S$ ; soit $g\in A_S$ et $g = g_1g_2^{-1}$ son \'ecriture normale \`a droite; alors les assertions suivantes sont \'equivalentes:\\(a) $gt = sg$;\\ (b) $\left \{ \begin{array}{l} tg_2 = g_2u, \\ sg_1 = g_1u \end{array}\right.$ pour un certain $u\in S$.   
\end{Lemme}
\begin{preuve} Il est clair que $(b)\Rightarrow (a)$. Soit $g\in A_S$ et $g = g_1g_2^{-1}$ son \'ecriture normale \`a droite. Supposons que $sg = gt$. On a alors $(sg_1)(tg_2)^{-1} = sg_1g_2^{-1}t^{-1} = g =  g_1g_2^{-1}$  et par le lemme \ref{lemclassi}(ii), il existe $u\in A_S^+$ tel que $sg_1 = g_1u$ et $tg_2 = g_2u$. On alors en particulier $\ell(u) = \ell(s) = 1$ et donc $u\in S$.\end{preuve}
\begin{Lemme}\label{prodrub} Soit $(A_S,S)$ un syst\`eme d'Artin-Tits et $X,Y,Z \subset S$. Soit $g$ est un $Y$-ruban-$X$ positif de $A_S^+$.\\(i) $g$ est r\'eduit-$X$ si et seulement si pour tout $s\in X$, $g$ est une cha\^\i ne-$s$. De plus dans ce cas, $g$ est une $t$-cha\^\i ne-$s$ avec $t\in Y$ et est aussi $Y$-r\'eduit.\\
(ii) si  $g$ est r\'eduit-$X$ et si $h$ est un $Z$-ruban-$Y$ positif r\'eduit-$Y$ alors $hg$ est un $Z$-ruban-$X$ positif r\'eduit-$X$.
\end{Lemme}
\begin{preuve} (i) Par d\'efinition d'un ruban, dire que $g$ est r\'eduit-$X$ est \'equivalent au fait que $g$ est pour tout $s\in X$ une $t$-cha\^\i ne-$s$ avec $t\in Y$. D'autre part, si $g$ n'est pas $Y$-r\'eduit, il existe $s\in Y$ et $g_1\in A^+_S$ tel que $g = sg_1$. Puisqu'il existe $t\in X$ tel que $sg = gt$, on obtient par simplification que $g = sg_1 = g_1t$ ; ce qui contredit le fait que $g$ est r\'eduit-$X$. Donc $g$ est $Y$-r\'eduit.\\ (2) Il en est de m\^eme pour $h$ et $hg$ car il est clair que $hg$ est un $Z$-ruban-$X$ positif. On conclut en remarquant que le produit d'une $u$-cha\^\i ne-$t$ par une $t$-cha\^\i ne-$s$ est une $u$-cha\^\i ne-$s$. \end{preuve}
\begin{Prop}
Soit $(A_S,S)$ un syst\`eme d'Artin-Tits de type sph\'erique et $X$ une partie de $S$. on pose $H_X = \{g = g_1g_2^{-1}\in A_S|$ $g_1$ et $g_2$ sont deux $X$-rubans-$X$ positifs et $X$-r\'eduits $\}$. Alors $H_X$ est un sous-groupe de $A_S$. 
\end{Prop}
\begin{preuve} On a $1\in H_X$ ; si $g\in H_X$, alors $g^{-1}$ est aussi dans $H_X$. Soit $g,z$ deux \'el\'ements de $H_X$. Ils s'\'ecrivent $g = g_1g_2^{-1}$ et $z= z_1z_2^{-1}$ avec $g_1,g_2,z_1,z_2$ quatre $X$-rubans-$X$ positifs et $X$-r\'eduits. On a $gz = g_1g_2^{-1}z_1z_2^{-1}$. De plus  $g_2^{-1}z_1 = v_1v_2^{-1}$ avec $v_1\land_\succ v_2 = 1$. D'apr\`es le lemme \ref{corZen}, $v_1$ et $v_2$ sont deux $X$-rubans-$X'$ positifs pour $X'$ une partie de $S$. En outre, $v_1$ et $v_2$ sont $X$-r\'eduits : supposons que $v_1$ ne soit pas $X$-r\'eduit ; alors il n'est pas non plus r\'eduit-$X'$.  Soit $s\in X'$ tel que $v_1\succ s$.  puisque $g_2v_1 = z_1v_2$ on a aussi $z_1v_2\succ s$. Comme $v_1\land_\succ v_2 = 1$, $s$ ne divise pas $v_2$ \`a droite et $v_2$ est une $t$-cha\^\i ne-$s$ avec $t\in X$ par le lemme \ref{prodrub}(i) ; de plus $z_1$ est r\'eduit-$X$, donc  c'est une cha\^\i ne-$t$ par le m\^eme lemme. Ceci implique que $z_1v_2$  est une cha\^\i ne-$s$, ce qui est contradictoire avec le fait qu'il est divisible \`a droite par $s$ ({\it cf.} proposition \ref{Prprcha2}). Donc $v_1$ est $X$-r\'eduit et par sym\'etrie, $v_2$ aussi. Finalement $gz = (g_1v_1)(z_2v_2)^{-1}$ avec $g_1v_1$ et $z_2v_2$ deux $X$-rubans-$X'$ positifs qui sont $X$-r\'eduits et $gz = R_1R_2^{-1}$ avec $R_1 = g_1v_1\overline{g_1v_1}$ et $R_2 = z_2v_2\overline{g_1v_1}$ o\`u $g\to \overline{g}$ d\'esigne l'unique anti-automorphisme de mono\"\i de qui fixe $S$. Enfin, $R_1$ et $R_2$ sont r\'eduits-$X$ par le lemme \ref{prodrub}(ii).\end{preuve}    

\begin{Prop}\label{Prcox2}

Soit $(A_S,S)$ un syst\`eme d'Artin-Tits de type sph\'erique et $X$ une partie de $S$. Alors $N_{A_S}(A_X) = H_X\ltimes A_X$.\end{Prop}

\begin{preuve} Puisque $A_S$ est de type sph\'erique, on  a $N_{A_S}(A_X) = A_X\cdot QZ_{A_S}(A_X)$ (\textit{cf.} th\'eor\`eme \ref{THdeb2}) ; d'autre part, on a $A_X\cap QZ_{A_S}(A_X) = QZ_{A_X}(A_X)$.\\ Il est clair que $QZ_{A_X}(A_X)\cap H_X = \{1\}$. D'autre part, tout \'el\'ement $g$ de $QZ_{A_S}(A_X)$ s'\'ecrit $g = g_1g_2^{-1}$ avec $g_1$ et $g_2$ deux $X$-rubans-$X'$ positifs premiers entre eux pour $\succ$. Mais par le lemme \ref{lemsph1} on peut \'ecrire $g_1 = \Delta_X^{n}g'_1$ et $g_2 = \Delta_X^mg'_2$ avec $n,m\in \mathbb{N}$ et $g'_1,g'_2$ deux $X$-rubans-$X'$ positifs qui sont $X$-r\'eduits. D'o\`u $g = \Delta_X^{n-m}g'_1{g'}_2^{-1}$ et $ g'_1{g'}_2^{-1}\in H_X$. Donc $QZ_{A_S}(A_X) = QZ_{A_X}(A_X)\rtimes H_X$.\end{preuve}

\begin{Lemme}\label{lemsurjHG} Soit $(A_S,S)$ un syst\`eme d'Artin-Tits de type sph\'erique et $X$ une partie de $S$. On note $p : A_S\to W_S$ le morphisme canonique. On a $p(H_X) = G_X$ o\`u $G_X = \{w\in W_S ; w\Pi_X = \Pi_X\}$. \end{Lemme}
\begin{preuve}
Soit $g\in Conj^+(S;X,X)$ r\'eduit-$X$. Par le th\'eor\`eme
 \ref{cordeb1}, on a $g = g_n\cdots g_1$ o\`u $g_i$ est un
 $X_i$-ruban-$X_{i-1}$\ positif \'el\'ementaire, $X_0 = X = X_n$. Notons
 $s_i,t_i\in S $ tels que $X_{i-1}\cup\{s_i\} = X_i\cup\{t_i\}$. Comme
 $g$ est r\'eduit-$X$, $g_i\not= \Delta_{X_i}$ donc $g_i = \Delta_{X_{i-1}(s_i)}\Delta_{X_{i-1}(s_i)-\{s_i\}}^{-1} = \Delta_{X_i(t_i)-\{t_i\}}^{-1}\Delta_{X_i(t_i)}$. Donc $p(g_i) = \nu(X_i,t_i)$ et $p(g) = \nu(X_n,t_n)\cdots \nu(X_1,t_1)$; donc $p(g)\in G_X$ (cette d\'ecomposition de $p(g)$ ne v\'erifie pas n\'ecessairement $\ell(p(g)) =\sum_{i=1}^{n}{\ell(\nu(X_i,t_i))}$~).  La r\'eciproque est identique \`a la preuve du th\'eor\`eme \ref{cordeb1} dans le cas o\`u $g$ est r\'eduit-$X$: soit $w\in G_X$ ; alors il existe une suite $s_0,\ldots,s_n$ d'\'el\'ements de $S$ et une suite $X_0 =X,X_1,\ldots,X_{n-1},X_n = X$ de parties de $S$ telles que $w = \nu(X_0,s_0)\cdots \nu(X_n,s_n)$ o\`u pour $i\in\{0,\ldots,n-1\}$, le but de $\nu(X_i,s_i)$ est $X_{i+1}$ et telles que $l(w) = \sum_{i=0}^{n-1}{\ell(\nu(X_i,s_i))}$. On a donc $\pi(\nu(X_i,s_i)) = \Delta_{X_i(s_i)-\{s_i\}}^{-1}\Delta_{X_i(s_i)} = \Delta_{X_{i+1}(t_i)}\Delta_{X_{i+1}(t_i)-\{t_i\}}^{-1}$ avec $t_i\in S$ tel que $X_{i+1}\cup\{t_i\} = X_i\cup\{s_i\}$. Donc  $\pi(w) \in H_X$. Il est r\'eduit-$X$ puisque $w$ l'est. Enfin $p(\pi(w)) = w$ par d\'efinition. 
\end{preuve}
\begin{preuve}[du th\'eor\`eme \ref{THcox}] Par la proposition \ref{Prcox2} on a l'isomorphisme $N_{A_S}(A_X)/(Z_{A_S}(A_X)\cdot A_X)\simeq H_X/(Z_{A_S}(A_X)\cap H_X),$ et par la proposition \ref{propkracox}, on a aussi $N_{W_S}(W_X)/(Z_{W_S}(W_X)\cdot W_X)\simeq G_X/(Z_{W_S}(W_X)\cap G_X).$ Par le lemme \ref{lemsurjHG}, $p(H_X) = G_X$ donc le morphisme $\widetilde{p}: H_X/(Z_{A_S}(A_X)\cap H_X)\to G_X/(Z_{W_S}(W_X)\cap G_X)$ induit par $p$ est aussi surjectif. Reste \`a voir qu'il est injectif. Mais, si $g\in H_X$ on a pour $s,t$ dans $X$ l'\'equivalence suivante \begin{equation} sg = gt \iff sp(g) = p(g)t\label{equation5}\end{equation} le sens $\Leftarrow$ vient du fait que pour $g$ dans $H_X$ et $s$ dans $X$, on a par d\'efinition $sg = gt'$ pour un certain $t'$ de $X$; par projection dans $W_S$, on trouve $sp(g) = p(g)t' = p(g)t$ et donc $t' = t$. Soit maintenant $\widetilde{g}$ dans $H_X/(Z_{A_S}(A_X)\cap H_X)$ de repr\'esentant $g$ dans $H_X$ et tel que $\widetilde{p}(\widetilde{g}) = 1$. Par d\'efinition, $\widetilde{p}(\widetilde{g}) = 1 \iff p(g)\in Z_{W_S}(W_X)\cap G_X$. Par l'\'equivalence (\ref{equation5}), cela entra\^\i ne $g\in Z_{A_S}(A_X)\cap H_X$ et donc  $\widetilde{g} = 1$. D'o\`u l'injectivit\'e de $\widetilde{p}$.\end{preuve}

\begin{tabular}{l} Eddy Godelle \\LAMFA CNRS 2270\\Universit\'e de Picardie-Jules
 Verne\\Facult\'e de Math\'ematiques et d'Informatique\\33 rue Saint-Leu
 80000 Amiens\\ France\\ eddy.godelle@u-picardie.fr \end{tabular}
\end{document}